\documentclass[12pt]{article}
\usepackage{amsmath}
\usepackage{amsfonts}
\usepackage{amssymb}
\usepackage[dvips]{graphicx}

\newif{\ifcomentarios}
\comentariosfalse
\newtheorem{theorem}{Theorem}

\newtheorem{lemma}[theorem]{Lemma}

\newtheorem{proposition}[theorem]{Proposition}
\newtheorem{remark}[theorem]{Remark}

\newcommand{\doint}{\displaystyle\oint}
\newcommand{\func}{\text}
\begin{document}

\author{\textbf{William R. P. Conti}\thanks{%
Address: Instituto de F\'{\i}sica, Universidade de S\~{a}o Paulo, Caixa
Postal 66318, 05314-970 S\~{a}o Paulo, SP, Brasil. Supported by FAPESP under
grant $\#07/59739-4$. E-mail: \textit{\ }\texttt{wrpconti@if.usp.br}}~ and \
~\textbf{Domingos H. U. Marchetti}\thanks{%
Present Address: Mathematics Department, \ The University of British
Columbia, Vancouver, \ BC, \ Canada \ V6T 1Z2. Email: \texttt{%
marchett@math.ubc.ca}} }
\title{Singular perturbation of nonlinear systems with regular singularity}
\date{}
\maketitle

\begin{abstract}
We extend Balser-Kostov method of studying summability properties of a
singularly perturbed inhomogeneous linear system with regular singularity at
origin to nonlinear systems of the form%
\begin{equation*}
\varepsilon zf^{\prime }=F(\varepsilon ,z,f)
\end{equation*}%
with $F$ a $\mathbb{C}^{\nu }$--valued function, holomorphic in a polydisc $%
\bar{D}_{\rho }\times \bar{D}_{\rho }\times \bar{D}_{\rho }^{\nu }$. We show
that its unique formal solution in power series of $\varepsilon $, whose
coefficients are holomorphic functions of $z$, is $1$--summable under a
Siegal--type condition on the eigenvalues of $F_{f}(0,0,0)$. The estimates
employed resemble the ones used in KAM theorem. A simple Lemma is developed
to tame convolutions that appears in the power series expansion of nonlinear
equations.

\medskip 


\end{abstract}

\section{Introduction}

\setcounter{equation}{0} \setcounter{theorem}{0}

Balser and Kostov\cite{BK} have studied singularly perturbed linear system
with regular singularity at $z=0$ of the form%
\begin{equation}
\varepsilon zf^{\prime }=Af-b  \label{syst}
\end{equation}%
$f^{\prime }$ means derivative of $f$ w.r.t. $z$; $A=A(\varepsilon ,z)$ and $%
b=b(\varepsilon ,z)$ are, respectively, a $\nu \times \nu $ matrix and a $%
\nu $--vector whose entries are holomorphic in the polydisc $D_{R}\times
D_{R}$, $R>0$.\footnote{%
Here, $D_{\rho }(z_{0})=\left\{ z\in \mathbb{C}:\left\vert
z-z_{0}\right\vert <\rho \right\} $ denotes an open disc of radius $\rho >0$%
, centered at $z_{0}$, $\bar{D}_{\rho }(z_{0})$ denotes its closure and $%
D_{\rho }=D_{\rho }(0)$.} $A$ is, in addition, such that $A(0,0)^{-1}$
exists. For such a system, there exists a unique formal solution in the ring 
$\mathcal{O}(r)[[\varepsilon ]]_{1}$ of formal power series 
\begin{equation}
\hat{f}(\varepsilon ,z)=\sum_{i=0}^{\infty }{a_{i}(z)\,\varepsilon ^{i}}
\label{fseries}
\end{equation}%
in $\varepsilon $ with coefficients ${a_{i}}(z)$ in the ring $\mathcal{O}(r)$
of holomorphic functions on $D_{r}$, continuous in its closure, satisfying 
\begin{equation}
\max_{\left\vert z-z_{0}\right\vert \leq r}\left\vert {a_{i}(z)}\right\vert
\leq C\mu ^{i}i!\ ,\qquad i=0,1,2,\ldots  \label{a_i}
\end{equation}%
for some positive constants $C$, $\mu $ and $0<r<R$. The authors have shown
(see Theorem 1 and 2 of \cite{BK}) that $\hat{f}(\varepsilon ,z)$ is the $1$%
--Gevrey asymptotic expansion of a holomorphic function $f(\varepsilon ,z)$
in $S(\theta ,\gamma ;E)\times D_{r}$, as $\varepsilon $ tends to $0$, if
the closed sector $\bar{S}(\theta ,\gamma ;E)$ does not contain any ray on
the direction of the eigenvalues $\lambda _{j}$ of $A(0,0)$: 
\begin{equation}
\left\vert \arg \lambda _{j}-\theta \right\vert >\gamma /2~\ ,\qquad
j=1,\ldots ,n~.  \label{lambda}
\end{equation}%
The formal series $\hat{f}(\varepsilon ,z)$ is thus $1$--summable in the
direction $\theta $ provided the eigenvalues of $A(0,0)$ satisfy a
Siegel--type condition, i.e. the $\lambda _{j}$ satisfy (\ref{lambda}) for
some $\gamma \geq \pi $.

A nonlinear version of (\ref{syst}) appears as follows. Let $f(\varepsilon
,z)$ be the unique extension in $S(0,\gamma ;E)\times D_{r}$, with $%
\varepsilon =2/N$, of the meromorphic function 
\begin{equation}
\phi _{\varepsilon }(z)=\frac{i}{2\sqrt{z}}\frac{J_{N/2}(i\sqrt{z}N)}{%
J_{N/2-1}(i\sqrt{z}N)}  \label{phie}
\end{equation}%
where $J_{\kappa }(x)$ is the Bessel function of order $\kappa $. This
function is related with the Fourier-Stieltjes transform $\hat{\sigma}%
^{N}(x) $ of a uniform measure $\sigma ^{N}$ on the $N$--dimensional sphere
of radius $\sqrt{N}$ and we refer to \cite{MC} and \cite{MCG} for the
motivations for its study. The $N$ dependence in the argument is chosen in
such way that $\phi _{\varepsilon }(z)$ attains, as $\varepsilon $ goes to $%
0 $, a limit function 
\begin{equation}
\phi _{0}(z)=\frac{-1}{1+\sqrt{1+4z}}~  \label{phi0}
\end{equation}%
(see Proposition 2.1 of \cite{MCG}). $\phi _{\varepsilon }$ satisfies an
ordinary (Riccati) differential equation 
\begin{equation}
\varepsilon z\phi _{\varepsilon }^{\prime }+\phi _{\varepsilon }-2z\phi
_{\varepsilon }^{2}+\frac{1}{2}=0~  \label{phieq}
\end{equation}%
which, despite of being nonlinear, can be dealt by Balser--Kostov's method.
It has been shown by the present authors (see Lemmas 3.2, 3.3 and 3.4 of 
\cite{MC}) \textbf{(a)} existence of a unique formal solution $\hat{\phi}%
_{\varepsilon }(z)$ in the form of (\ref{fseries}), satisfying (\ref{a_i}); 
\textbf{(b)} $\hat{\phi}_{\varepsilon }(z)$ is the $1$--Gevrey asymptotic
expansion of the holomorphic solution $f(\varepsilon ,z)$ of (\ref{phieq})
in $S(0,\gamma ;E)\times D_{r}$, as $\varepsilon $ goes to $0$ in $%
S(0,\gamma ;E)$; \textbf{(c)} choosing the sector $S(\theta ,\gamma ;E)$ of
opening angle $\gamma >\pi $ away from the negative real axis, $\hat{\phi}%
_{\varepsilon }(z)$ is, in addition, $1$-- summable in $\theta $ direction
and its sum is equal to $f(\varepsilon ,z)$.

In the present article statements \textbf{(a)}--\textbf{(c)}, together with
the $1$-summability, will be extended for more general ordinary differential
equations of the form%
\begin{equation}
\varepsilon zf^{\prime }=F(\varepsilon ,z,f)~,  \label{F}
\end{equation}%
with $f=(f^{1},\ldots ,f^{\nu })$ and $F=(F^{1},\ldots ,F^{\nu })$ $\nu $%
--vector functions, $F^{i}$ holomorphic in a polydisc, say $\bar{D}_{\rho
}\times \bar{D}_{\rho _{1}}\times \bar{D}_{\rho }^{\nu }$, for some $\rho
_{1}>\rho >0$. As in (\cite{BK}), the $\nu \times \nu $ matrix $%
A_{0,1}(0)=F_{f}(0,0,0)$ is assumed to be invertible, a condition that makes
(\ref{F}) to possess a regular singularity at $z=0$, and every eigenvalue of 
$A_{0,1}(0)$ satisfies condition (\ref{lambda}). Equation (\ref{phieq}) is
of the form (\ref{F}) with $\nu =1$ and\footnote{%
Statements \textbf{(a)}--\textbf{(c)} hold with $1/2$ in (\ref{phieq})
replaced by $\beta (\varepsilon )/2$ for any $1$--summable $\beta
(\varepsilon )=\displaystyle\sum_{n\geq 0}\beta _{n}\varepsilon ^{n}$ formal
series in $\theta $ direction. In this case the limit function (\ref{phi0})
is replaced by $\phi _{0}=-\beta _{0}/(1+\sqrt{1+4\beta _{0}z})$.} 
\begin{equation}
F(\varepsilon ,z,f)=-\frac{\beta (\varepsilon )}{2}-f+2zf^{2}\   \label{FF}
\end{equation}

Balser--Kostov summability proof in \cite{BK} of the formal series $\hat{f}$
solution does not follow the usual route: the (formal) Borel transform $%
\mathcal{\hat{B}}\hat{f}$ of $\hat{f}$ is analytically continued along some
sector of infinite radius (see e.g. \cite{Ba}). Their proof establishes
instead Grevrey asymptotic expansion directly from the equation (\ref{syst}%
), making resource of an auxiliary Lemma regarding an infinite system of
linear equation of the same type. Although (\ref{F}) is nonlinear, the
system of infinitely many equations obtained by taking derivatives of (\ref%
{F}) with respect to $\varepsilon $ is linear, indeed of the type stated in
Lemma 3 of \cite{BK}, and Balser--Kostov's method carries over to equation
of the form (\ref{F}).

The layout of this paper is as follows. In Section \ref{PS} we prove
existence of a unique solution of (\ref{F}) in power series of $z$. In
Section \ref{FPS} we show that the formal power series solution of (\ref{F})
is Gevrey of order $1$. In Section \ref{GA} Gevrey asymptotics are
established. Our main result, the $1$--summability of the formal solution of
(\ref{F}), is stated in Section \ref{S} and proved using Propositions \ref%
{P1}-\ref{P3} of the previous sections. The main ingredient (Lemma \ref{key}%
), is employed to tame arbitrarily large number of convolutions arised in
the expansion of $F$ in power series of $f$.

\section{Power series in $z$\label{PS}}

\setcounter{equation}{0} \setcounter{theorem}{0}

Under the hypothesis on $F$, the series

\begin{equation}
F(\varepsilon ,z,f)=\sum_{\substack{ n,m=0:  \\ n+m\neq 0}}^{\infty
}A_{n,m}(\varepsilon )z^{n}f^{m}  \label{Fab}
\end{equation}%
\newline
converges (in norm) absolutely in $\bar{D}_{\rho _{1}}\times \bar{D}_{\rho
}^{\nu }$, uniformly in $\varepsilon \in \bar{D}_{\rho }$, with the
coefficients $A_{n,m}(\varepsilon )$ regarded as a multilinear operator, 
\begin{equation*}
f^{m}\in \underset{m\ \mathrm{copies}}{\underbrace{\mathbb{C}^{\nu }\times
\cdots \times \mathbb{C}^{\nu }}}\longmapsto A_{n,m}(\varepsilon )f^{m}\in 
\mathbb{C}^{\nu }
\end{equation*}%
\begin{equation}
\left( A_{n,m}(\varepsilon )f^{m}\right) ^{i}=\sum_{i_{1},\ldots
,i_{m}=1}^{\nu }A_{n,m}^{i,i_{1},\ldots ,i_{m}}(\varepsilon )f^{i_{1}}\cdots
f^{i_{m}}~,  \label{Anm}
\end{equation}%
endowed with an operator norm induced by the Euclidean space $\mathbb{C}%
^{\nu }$:%
\begin{equation*}
\left\Vert A_{n,m}(\varepsilon )\right\Vert =\sup_{(v_{1},\ldots ,v_{m})\in 
\mathbb{C}^{m\nu }}\frac{\left\Vert A_{n,m}(\varepsilon )v_{1}\cdots
v_{m}\right\Vert }{\left\Vert v_{1}\right\Vert \cdots \left\Vert
v_{m}\right\Vert }~,
\end{equation*}%
holomorphic in $D_{\rho }$ as a function of $\varepsilon $.

In (\ref{Anm}) and from now on, $f=\left( f^{1},\ldots ,f^{\nu }\right) $
denotes a $\nu $--vector with $i$--th component $f^{i}$ and Euclidean norm $%
\left\Vert f\right\Vert ^{2}=f\cdot f=\displaystyle\sum_{i=1}^{\nu }\bar{f}%
^{i}f^{i}$. Without loss of generality, we assume $A_{00}(\varepsilon
)\equiv 0$ and since the l.h.s. of (\ref{F}) vanishes for $z=0$, its
solution in power series reads 
\begin{equation}
f(\varepsilon ,z)=\sum_{k=1}^{\infty }f_{k}(\varepsilon )z^{k}~.
\label{ffseries}
\end{equation}%
(by hypothesis $f(\varepsilon ,0)\equiv 0$). For $F$ given by the example (%
\ref{FF}), $A_{0,0}(\varepsilon )=\beta (\varepsilon )/2$ does not vanishes
and we may replace $f$ and (\ref{FF}) by $\tilde{f}=f+\beta /2$ and%
\begin{equation*}
\tilde{F}=\left( 1-2\beta z\right) \tilde{f}+\frac{\beta ^{2}}{2}z+2z\tilde{f%
}^{2}
\end{equation*}%
which satisfy $\tilde{f}(\varepsilon ,0)=0$ and $\tilde{A}_{0,0}(\varepsilon
)=0$. The general case differs very little from this particular example.

Substituting the power series (\ref{ffseries}) into (\ref{Fab}) together
with (\ref{F}), we are led to a system of equations%
\begin{equation}
\left( \varepsilon jI-A_{0,1}(\varepsilon )\right) f_{j}=g_{j}(\varepsilon
;f_{1},\ldots ,f_{j-1})  \label{fk}
\end{equation}%
with $g_{j}=(g_{j}^{1},\ldots ,g_{j}^{\nu })$ given by%
\begin{equation*}
g_{1}^{i}(\varepsilon )=A_{1,0}^{i}(\varepsilon )
\end{equation*}%
and%
\begin{equation}
g_{j}^{i}(\varepsilon ;f_{1},\ldots ,f_{j-1})=\sum_{\substack{ n,m:  \\ %
2\leq n+m\leq j}}\sum_{i_{1},\ldots ,i_{m}=1}^{\nu }A_{n,m}^{i,i_{1},\ldots
,i_{m}}(\varepsilon )\left( f^{i_{1}}\ast \cdots \ast f^{i_{m}}\right) _{j-n}
\label{sum}
\end{equation}%
for $j\geq 2$; for any two sequences $\alpha =\left( \alpha _{k}\right)
_{k\geq 1}$ and $\beta =\left( \beta _{k}\right) _{k\geq 1}$, their
convolution product $\alpha \ast \beta =\left( (\alpha \ast \beta
)_{k}\right) _{k\geq 1}$ is a sequence defined by $\left( \alpha \ast \beta
\right) _{1}=0$ and 
\begin{equation}
\left( \alpha \ast \beta \right) _{k}=\sum_{l=1}^{k-1}\alpha _{l}\beta
_{k-l}~,\quad k\geq 2~.  \label{prod}
\end{equation}%
The restriction $n+m\leq j$ in (\ref{sum}) results from the fact that our
sequence $f^{i}=\left( f_{k}^{i}\right) _{k\geq 1}$ starts with $k=1$ and a
convolution involving $m$ sequences cannot have nonvanishing component $j-n$
if $j>n+m$.

Consequently, for any $k\in \mathbb{N}$ arbitrary, (\ref{fk}) for $1\leq
j\leq k$ forms a \textbf{closed} system of $\nu \cdot k$ equations,
involving $\nu \cdot k$ unknown functions which can be solved by iteration
starting from%
\begin{equation}
f_{1}(\varepsilon )=\left( \varepsilon I-A_{0,1}(\varepsilon )\right)
^{-1}A_{1,0}(\varepsilon )~.  \label{fg1}
\end{equation}%
If equation (\ref{fk}) for $1\leq j\leq k-1$ and $k\geq 2$ have been solved,
then%
\begin{equation}
f_{k}(\varepsilon )=\left( \varepsilon kI-A_{0,1}(\varepsilon )\right)
^{-1}g_{k}(\varepsilon ;f_{1},\ldots ,f_{k-1})\ .\qquad  \label{fgk}
\end{equation}%
Regarding the inverse matrix $\left( \varepsilon kI-A_{0,1}(\varepsilon
)\right) ^{-1}$, we have the following

\begin{lemma}[Lemma 1 of \protect\cite{BK}]
Suppose (\ref{lambda}) holds with $\theta =0$ and $\lambda _{j}$, $%
j=1,\ldots ,\nu $, eigenvalues of $A_{0,1}(0)=F_{f}(0,0,0)$. One can always
find $E>0$ such that, if $k\left\vert \varepsilon \right\vert \geq \dfrac{c}{%
c-1}\sup_{\left\vert \varepsilon \right\vert \leq E}\left\Vert
A_{0,1}(\varepsilon )\right\Vert $ for some $c>1$, the inverse matrix in (%
\ref{fgk}), given by 
\begin{equation*}
\left( \varepsilon kI-A_{0,1}(\varepsilon )\right) ^{-1}=\sum_{n=0}^{\infty }%
\frac{1}{\left( \varepsilon k\right) ^{n+1}}\left( A_{0,1}(\varepsilon
)\right) ^{n}
\end{equation*}%
is bounded and satisfies $\left\Vert \left( \varepsilon
kI-A_{0,1}(\varepsilon )\right) ^{-1}\right\Vert \leq c$, uniformly in $%
S\left( 0,\gamma ;E\right) $. If $k\left\vert \varepsilon \right\vert <%
\dfrac{c}{c-1}\sup_{\left\vert \varepsilon \right\vert \leq E}\left\Vert
A_{0,1}(\varepsilon )\right\Vert $, let $\lambda _{j}(\varepsilon )$, $%
j=1,\ldots ,n$, the eigenvalues of $A_{0,1}(\varepsilon )$, be so that their
distances from every ray $\eta =re^{i\tau }$ intercepting $S\left( 0,\gamma
;E\right) $ are bounded from below by a constant $a>0$:%
\begin{equation}
a=\inf \left\{ \left\vert \lambda _{j}(\varepsilon )-re^{i\tau }\right\vert
:0\leq r<\infty \ ,\ \left\vert \tau \right\vert \leq \gamma ~,~j=1,\ldots
,n~\text{and }\varepsilon \in S\left( 0,\gamma ;E\right) \right\} ~.
\label{a}
\end{equation}%
Then, 
\begin{equation*}
\left\vert \det \left( \varepsilon kI-A_{0,1}(\varepsilon )\right)
\right\vert =\prod_{j=1}^{\nu }\left\vert \varepsilon k-\lambda
_{j}(\varepsilon )\right\vert \geq a^{\nu }>0
\end{equation*}%
together with the formula $A^{-1}=\func{Adj}(A)/\det A$ for inverse of a
matrix $A$, where $\func{Adj}(A)$ is the transposed of the cofactors matrix
of $A$, (see e.g. \cite{La}) and with the boundedness in $S\left( 0,\gamma
;E\right) $ of all cofactors of $A_{0,1}(\varepsilon )$, give 
\begin{equation}
\left\Vert \left( \varepsilon kI-A_{0,1}(\varepsilon )\right)
^{-1}\right\Vert \leq c~,  \label{c}
\end{equation}%
uniformly in $S\left( 0,\gamma ;E\right) $ for every $k\in \mathbb{N}$.
\end{lemma}

\begin{proposition}
\label{P1} Let $F$ be given by (\ref{Fab}) with the eigenvalues of $%
A_{0,1}(0)$ obeying hypothesis (\ref{lambda}). There exist $\gamma $, $E$
and $\sigma $ such that (\ref{F}) has a solution $f(\varepsilon ,z)$
holomorphic in $S(0,\gamma ;E)\times {}D_{\sigma }$. The solution $%
f(\varepsilon ,z)$ converges, as $\varepsilon \rightarrow 0$ in the sector $%
S(0,\gamma ;E)$, to the unique solution $f^{\ast }(z)$ of $F(0,z,f)=0$ in $%
D_{\sigma }$ satisfying $f(0)=0$.
\end{proposition}

\noindent \textbf{Proof} Since (\ref{ffseries}) solves (\ref{F}), its
coefficients $f_{k}(\varepsilon )$ satisfy the formal relations (\ref{fk})
whose solution depends on the existence of inverse matrix $\left(
\varepsilon kI-A_{0,1}(\varepsilon )\right) ^{-1}$ for every $k\in \mathbb{N}
$ and $\varepsilon \in S\left( 0,\gamma ,E\right) $. Assuming (\ref{lambda})
holds for every eigenvalues of $A_{0,1}(0)$, let $\gamma $ and $E$ be such
that (\ref{a}), and consequently (\ref{c}), holds. Hence, $f_{k}(\varepsilon
)$ given by (\ref{fgk}) is bounded uniformly in $S\left( 0,\gamma ;E\right) $%
, uniquely defined for every $k\in \mathbb{N}$ and, in view of these,
holomorphic in $S\left( 0,\gamma ;E\right) $.

Let $\phi _{l}$ and $\alpha _{n,m}$ be the supremum in $S(0,\gamma ;E)$ of $%
\left\Vert f_{l}(\varepsilon )\right\Vert $ and $\left\Vert
A_{n,m}(\varepsilon )\right\Vert $, respectively:%
\begin{eqnarray}
\phi _{l} &=&\sup_{\varepsilon \in S(0,\gamma ;E)}\left\Vert
f_{l}(\varepsilon )\right\Vert  \notag \\
\alpha _{n,m} &=&\sup_{\varepsilon \in S(0,\gamma ;E)}\left\Vert
A_{n,m}(\varepsilon )\right\Vert ~.  \label{alphanm}
\end{eqnarray}%
By Cauchy formula 
\begin{equation*}
\frac{1}{n!m!}F^{(0,n,m)}(\varepsilon ,0,0)\left( \frac{f}{\left\Vert
f\right\Vert }\right) ^{m}=\frac{1}{(2\pi i)^{2}}\doint \doint \frac{%
F(\varepsilon ,\zeta ,\phi f/\left\Vert f\right\Vert )}{\zeta ^{n+1}\phi
^{m+1}}d\zeta d\phi
\end{equation*}%
and there exists $C<\infty $ ($=\sup_{\bar{S}(0,\gamma ;E)\times {}\bar{D}%
_{\rho _{1}}\times \bar{D}_{\rho }^{\nu }}\left\Vert F(\varepsilon
,z,f)\right\Vert $, $E\leq \rho $) such that%
\begin{equation}
\alpha _{n,m}\leq \frac{C}{\rho _{1}^{n}\rho ^{m}}~.  \label{alphab}
\end{equation}%
Now, we prove that the majorant series $\displaystyle\sum_{l=1}^{\infty
}\phi _{l}\sigma ^{l}$ converges and is bounded by $\rho $ for some $%
0<\sigma <\rho $. For this, the following lemma will be stated more
generally than it is needed for this section.

\begin{lemma}
\label{key}\vspace{1mm}Let $\lambda \geq 0$ be given and let $A=\left( 1+\pi
^{2}/3\right) ^{-1}/2=0.1165536\ldots $. Consider the sequence $\left(
C_{l}\right) _{l=0}^{\infty }$ with $C_{0}=A$ or $C_{0}=0$ and \footnote{%
The sequence $\left( C_{l}\right) $ of this Section has $C_{0}=0$. Lemma \ref%
{key} has been stated with $C_{0}=A$ to be used elsewhere in another Section.%
}%
\begin{equation*}
C_{l}=\frac{Al!^{\lambda }}{l^{2}}\ ,\hspace{10mm}\forall \;l\geq 1~.
\end{equation*}%
Then%
\begin{equation}
\sum_{l=0}^{m}{C_{l}\,C_{m-l}}\leq {}C_{m}  \label{conv}
\end{equation}%
\noindent holds for every $m\geq 0$.
\end{lemma}

\noindent \textbf{Proof }Since $\dbinom{m}{l}\geq 1$,%
\begin{equation}
\frac{1}{l}+\frac{1}{m-l}=\frac{m}{l(m-l)}  \notag
\end{equation}%
and $0\leq (a-b)^{2}=2(a^{2}+b^{2})-(a+b)^{2}$ holds for any real numbers $a$
and $b$, we have%
\begin{eqnarray*}
\frac{1}{C_{m}}\sum_{l=0}^{m}{C_{l}\,C_{m-l}} &\leq &A\left(
2+\sum_{l=1}^{m-1}{\frac{m^{2}}{l^{2}(m-l)^{2}}}\right) \\
&\leq &2A\left( 1+\sum_{l=1}^{m-1}{\left( \frac{1}{l^{2}}+\frac{1}{(m-l)^{2}}%
\right) }\right) \\
&\leq &2A\left( 1+\frac{\pi ^{2}}{3}\right) =1~.
\end{eqnarray*}

\hfill $\Box $

It thus follows from (\ref{conv}) with $C_{0}=0$ that 
\begin{equation}
\sum_{\underset{l_{1}+\cdots +l_{k}=m}{l_{1},\ldots ,l_{k}\geq 1:}%
}C_{l_{1}}\cdots C_{l_{k}}\leq {C_{m}}\;  \label{convk}
\end{equation}%
holds for any $1\leq {}k\leq {}m$.

Let us assume that (\ref{alphab}) can be written as (see Remark \ref{R1}) 
\begin{equation}
\alpha _{n,m}\leq \frac{\alpha }{c}C_{n}\frac{1}{\rho ^{n+m}}
\label{alpha_nm}
\end{equation}%
and suppose%
\begin{equation*}
\phi _{l}\leq \alpha {}C_{l}\frac{1}{\kappa ^{l}}
\end{equation*}%
holds for $l\geq 1$, with $\left( C_{l}\right) _{l\geq 1}$ the sequence in
Lemma \ref{key} with $\lambda =0$, for some $\alpha $ and $\kappa <\rho $.
Hence, by (\ref{fg1}) together with (\ref{c}) and (\ref{alphanm}), we have%
\begin{equation}
{}\phi _{1}\leq c\alpha _{1,0}\leq \alpha {}A\frac{1}{\kappa }  \label{phi1}
\end{equation}%
and, by (\ref{fgk}) and (\ref{sum}) together with (\ref{c}),%
\begin{equation*}
\left\Vert f_{k}(\varepsilon )\right\Vert \leq c\sum_{\substack{ n,m:  \\ %
2\leq n+m\leq k}}\left\Vert A_{n,m}(\varepsilon )\right\Vert \left(
\left\Vert f(\varepsilon )\right\Vert \ast \cdots \ast \left\Vert
f(\varepsilon )\right\Vert \right) _{k-n}
\end{equation*}%
with $\left\Vert f(\varepsilon )\right\Vert $ denoting the sequence $\left(
\left\Vert f_{j}(\varepsilon )\right\Vert \right) _{j\geq 1}$. Taking the
supremum over $\varepsilon \in S(0,\gamma ;E)$ in both sides together with (%
\ref{alphanm}), (\ref{alphab}) and (\ref{convk}), 
\begin{eqnarray*}
\phi _{k} &\leq &c\left( \alpha _{k,0}+\sum_{\substack{ n,m\geq 1:  \\ 2\leq
n+m\leq k}}\alpha _{n,m}\left( \underset{m}{\underbrace{\phi \ast \cdots
\ast \phi }}\right) _{k-n}\right) \\
&\leq &\alpha C_{k}\frac{1}{\rho ^{k}}+\alpha \frac{1}{\kappa ^{k}}%
\sum_{1\leq n\leq k-1}C_{n}\left( \frac{\kappa }{\rho }\right)
^{n}C_{k-n}\sum_{1\leq m\leq k-n}\left( \frac{\alpha }{\rho }\right) ^{m} \\
&\leq &\alpha \left( \frac{\kappa ^{k}}{\rho ^{k}}+\frac{\alpha }{\rho
-\alpha }\right) C_{k}\frac{1}{\kappa ^{k}}\leq \alpha C_{k}\frac{1}{\kappa
^{k}}\ ,
\end{eqnarray*}%
holds for $k\geq 2$ provided $\alpha <\rho /2$ and 
\begin{equation}
\kappa =\rho \sqrt{1-\frac{\alpha }{\rho -\alpha }}  \label{kappa}
\end{equation}%
With $\alpha $ and $\kappa $ satisfying these conditions, we conclude%
\begin{equation}
\phi _{l}=\sup_{\varepsilon \in {}S(0,\gamma ;E)}\left\Vert
f_{l}(\varepsilon )\right\Vert \leq \alpha \frac{A}{l^{2}}\frac{1}{\kappa
^{l}}~,\hspace{10mm}\forall \,l\geq 1  \label{phil}
\end{equation}%
and $\left( f_{l}(\varepsilon )\,z^{l}\right) _{l\geq 1}$ is a sequence of
holomorphic functions, uniformly bounded in $S(0,\gamma ;E)\times
{}D_{\sigma }$ by $\left( {\phi }_{l}{\,\sigma ^{l}}\right) _{l\geq 1}$,
whose sum $f(\varepsilon ,z)=\displaystyle\sum_{l=1}^{\infty
}f_{l}(\varepsilon )\,z^{l}$ is bounded (in norm) by%
\begin{equation}
\sum_{l=1}^{\infty }{\phi }_{l}{\,\sigma ^{l}}=\frac{\alpha A\kappa }{\kappa
-\sigma }=\rho  \label{maj}
\end{equation}%
provided $\sigma <\kappa $ satisfies $\sigma =\kappa (\rho -\alpha A)/\rho
=(\rho -\alpha A)\sqrt{1-\alpha /(\rho -\alpha )}$, by (\ref{kappa}). Under
this choice of $\sigma $, $F(\varepsilon ,z,D_{\rho }^{\nu })\subset
D_{\sigma }^{\nu }$ uniformly in $S(0,\gamma ;E)\times {}D_{\sigma }$ and
the solution we have obtained by the formal expansion (\ref{fk}) and (\ref%
{sum}) acquires sense. The power series solution (\ref{ffseries}) of (\ref{F}%
) thus converges to a unique analytic function $f(\varepsilon ,z)$ in $%
S(0,\gamma ;E)\times {}D_{\sigma }$. The proof of uniqueness will be omitted.

From the uniform convergence of (\ref{ffseries}) we conclude that, for any
fixed $z\in {}D_{\sigma }$, the solution $f(\varepsilon ,z)$ tends to%
\begin{equation*}
f(0,z)=\lim_{S(0,\gamma ;E)\ni \varepsilon \rightarrow 0}\sum_{l=0}^{\infty }%
{f_{l}(\varepsilon )\,z^{l}}=\sum_{l=0}^{\infty }\lim_{S(0,\gamma ;E)\ni
\varepsilon \rightarrow 0}{f_{l}(\varepsilon )\,z^{l}}=f^{\ast }(z)
\end{equation*}%
where $f^{\ast }(z)$ is the unique solution of equation 
\begin{equation*}
F(0,z,f)=0
\end{equation*}%
for $f$, by the analytic implicit function theorem (see e.g. Section 2.3 of 
\cite{Be} or the next section, for an alternative solution). Note that the
solution $f^{\ast }(z)$ is regular at $z=0$ since, by (\ref{ffseries}), it
must satisfy $f(0)=0$ and this concludes the proof of Proposition \ref{P1}.

\hfill $\Box $

\section{Formal power series in $\protect\varepsilon $\label{FPS}}

\setcounter{equation}{0} \setcounter{theorem}{0}

As in (\ref{Fab}), the double series%
\begin{equation}
F(\varepsilon ,z,f)=\sum_{n,m}B_{n,m}(z){\varepsilon ^{n}}f^{m}  \label{FB}
\end{equation}%
\newline
converges (in norm) absolutely in $\bar{D}_{\rho }\times \bar{D}_{\rho
}^{\nu }$, uniformly in $z\in \bar{D}_{\rho _{1}}$, with the coefficients $%
B_{n,m}(\varepsilon )$ regarded as a multilinear operator $f^{m}\in \mathbb{C%
}^{m\nu }\longmapsto B_{n,m}(z)f^{m}\in \mathbb{C}^{\nu }$%
\begin{equation*}
\left( B_{n,m}(z)f^{m}\right) ^{i}=\sum_{i_{1},\ldots ,i_{m}=1}^{\nu
}B_{n,m}^{i,i_{1},\ldots ,i_{m}}(z)f^{i_{1}}\cdots f^{i_{m}}~.
\end{equation*}%
By consistency, $B_{00}(0)=0$ but $B_{00}(z)$ may not be identically zero.
Before we go through the power series in $\varepsilon $, we study the
solution $f^{\ast }(z)$ of%
\begin{equation}
0=\sum_{m=0}^{\infty }B_{0,m}(z)a_{0}^{m}(z)~,  \label{B0}
\end{equation}%
in power series of $z$:%
\begin{equation}
a_{0}(z)=\sum_{j=1}^{\infty }a_{0,j}z^{j}~.  \label{a0}
\end{equation}%
Note that $f^{\ast }(z)=a_{0}(z)$, by Proposition \ref{P1}, so $a_{0}(0)=0$.
Replacing (\ref{a0}) into (\ref{B0}), and taking into account%
\begin{equation*}
B_{0,m}(z)=\sum_{n=0}^{\infty }A_{n,m}(0)z^{n}
\end{equation*}%
equation (\ref{B0}) can be written as (omitting the argument $\varepsilon =0$
of $A_{n,m}(0)$, for simplicity)%
\begin{equation*}
0=A_{j,0}+\sum_{1\leq m\leq j}\left( A_{\cdot ,m}\ast a_{0}\ast \cdots \ast
a_{0}\right) _{j}~.
\end{equation*}%
For $j=1$,%
\begin{equation*}
A_{1,0}+A_{0,1}a_{0,1}=0\Longrightarrow a_{0,1}=-A_{0,1}^{-1}A_{1,0}~.
\end{equation*}%
Now, supposing $a_{0,1},\ldots ,a_{0,k-1}$ have already been determined,
then 
\begin{equation}
a_{0,k}=-A_{0,1}^{-1}\left(
A_{k,0}+\sum_{j=1}^{k-1}A_{j,1}a_{0,k-j}+\sum_{2\leq m\leq k}\left( A_{\cdot
,m}\ast a_{0}\ast \cdots \ast a_{0}\right) _{k}\right) ~.  \label{a0k}
\end{equation}%
If one takes the norm of (\ref{a0k}), together with $\left\Vert
A_{0,1}^{-1}\right\Vert \leq c$, (\ref{alphab}), (\ref{convk}) and (\ref%
{phil}), that holds also for $\varepsilon =0$,%
\begin{equation}
\left\Vert a_{0,k}\right\Vert \leq cC\left( \frac{1}{\rho _{1}^{k}}%
+\sum_{j=1}^{k-1}\left( \frac{1}{\rho _{1}}\right) ^{j}\sum_{1\leq m\leq
k-j}\left( \frac{1}{\rho }\right) ^{m}\left( \left\Vert a_{0,\cdot
}\right\Vert \ast \cdots \ast \left\Vert a_{0,\cdot }\right\Vert \right)
_{k-j}\right) \leq \alpha C_{k}\frac{1}{\kappa ^{k}}  \label{a0knorm}
\end{equation}%
provided we fix $\alpha $ and $\kappa $ as in the previous section, which is
consistent with the domain in which $f(\varepsilon ,z)$ is holomorphic. This
shows that $f^{\ast }(z)$ is holomorphic in $D_{\sigma }^{\nu }$ and proves
the existence of a unique solution of $F(0,z,f)=0$ in the same domain.

\begin{remark}
\label{R1}Regarding the radius of convergence of the power series of $%
f^{\ast }(z)$ one can estimate it a little better using the Cauchy majorant
method as in Section 3.2 of \cite{Be} (see also \cite{BK}, Section 1, for
the linear equation). Multiplying (\ref{a0knorm}) by $z^{k}$, summing over $%
k $ and replacing the inequality by equality, yields%
\begin{equation*}
\phi (z)=cC\frac{z/\rho _{1}}{1-z/\rho _{1}}\frac{1}{1-\phi (z)/\rho }
\end{equation*}%
for a majorant $\phi (z)$ of $f^{\ast }(z)$, whose solution%
\begin{equation*}
\phi (z)=\frac{\rho }{2}\left( 1-\sqrt{\frac{1-z/\sigma _{1}}{1-z/\rho _{1}}}%
\right)
\end{equation*}%
is holomorphic in a disc $D_{\sigma _{1}}$ with $\sigma _{1}=\rho _{1}\rho
/(\rho +4cC)<\rho _{1}$, proportional to $\rho _{1}$. In Section \ref{PS},
we have chosen $\rho _{1}$ so large that (\ref{alphab}) can be written as (%
\ref{alpha_nm}) and the radius of convergence $\sigma $, obtained applying
Lemma \ref{key} to convolutions, is proportional to $\rho $ instead (see
expression after (\ref{maj})). Despite of this loss, the method introduced
there is undeniably practical, more adaptable to diverse situations and, for
these reasons, we shall apply it here and in further sections.
\end{remark}

\begin{proposition}
\label{P2} Suppose the formal power series (\ref{fseries}) satisfies
equation (\ref{F}), formally, with $F=F(\varepsilon ,z,f)$ obeying the
hypotheses stated after (\ref{F}). Then, the coefficients $\left(
a_{i}(z)\right) _{i\geq 0}$ of (\ref{fseries}) are analytic functions of $z$
in the open disc $D_{\kappa }$ and there exist positive constants $C$ and $%
\mu $ such that%
\begin{equation}
\left\Vert a_{i}(z)\right\Vert \leq {}Ci!\mu ^{i}  \label{1gevrey}
\end{equation}%
holds for all $i\geq 0$ and $z\in \bar{D}_{\sigma }$, with $\sigma <\kappa
<\rho $. In other words, the formal power series is of Gevrey order $1$,
i.e., $\hat{f}(\varepsilon ,z)\in \mathcal{O}(\sigma )[[\varepsilon ]]_{1}$.
\end{proposition}

\noindent \textbf{Proof} Substituting the power series (\ref{fseries}) into (%
\ref{FB}), we are thus led to a system of equations%
\begin{equation*}
0=\sum_{m=0}^{\infty }B_{0,m}(z)a_{0}^{m}(z)~,
\end{equation*}%
which has already been solved for $a_{0}(z)$, and for $i\geq 1$%
\begin{equation}
za_{i-1}^{\prime }(z)=\sum_{m=1}^{\infty
}mB_{0,m}(z)a_{0}^{m-1}(z)a_{i}(z)+\sum_{n=1}^{i}\sum_{m=1}^{\infty
}B_{n,m}\left( \underset{m}{\underbrace{a(z)\ast \cdots \ast a(z)}}\right)
_{i-n}~.  \label{Ba}
\end{equation}%
We observe that the sum over $m$ has no limit as the sequence $a(z)=\left(
a_{k}(z)\right) _{k\geq 0}$ starts from $k=0$ and the convolution product of
any two sequences $\alpha =\left( \alpha _{k}\right) _{k\geq 0}$ and $\beta
=\left( \beta _{k}\right) _{k\geq 0}$ is now defined by 
\begin{equation}
\left( \alpha \ast \beta \right) _{k}=\sum_{l=0}^{k}\alpha _{l}\beta
_{k-l}~,\quad k\geq 0~.  \label{conv1}
\end{equation}

To isolate $a_{i}$, the largest index term in (\ref{Ba}), we have to show
that the matrix (recall $B_{0,1}(0)=A_{0,1}(0)$) 
\begin{eqnarray}
T_{0}(z) &=&B_{0,1}(z)+\sum_{m=2}^{\infty }mB_{0,m}(z)a_{0}^{m-1}  \notag \\
&=&A_{01}(0)\left( I+A_{01}(0)^{-1}\left(
B_{0,1}(z)-B_{0,1}(0)+\sum_{m=2}^{\infty }mB_{0,m}(z)a_{0}^{m-1}\right)
\right)  \label{T0}
\end{eqnarray}%
is invertible for every $z\in D_{\kappa }$ for some $\kappa \leq \rho $. For
this, we take $\kappa $ so small that%
\begin{equation*}
c\sup_{z\in D_{\kappa }(0)}\left( \left\Vert
B_{0,1}(z)-B_{0,1}(0)\right\Vert +\sum_{m=2}^{\infty }m\left\Vert
B_{0,m}(z)\right\Vert \left\Vert a_{0}\right\Vert ^{m-1}\right) \leq b<1
\end{equation*}%
and, consequently, $\left\Vert T_{0}(z)^{-1}\right\Vert \leq c/(1-b)$ holds
uniformly in $D_{\kappa }(0)$.

It follows from (\ref{Ba}) and (\ref{T0}) that 
\begin{equation}
a_{i}(z)=T_{0}(z)^{-1}\left( za_{i-1}^{\prime
}(z)-\sum_{n=1}^{i}\sum_{m=1}^{\infty }B_{n,m}\left( \underset{m}{%
\underbrace{a(z)\ast \cdots \ast a(z)}}\right) _{i-n}\right)  \label{ai}
\end{equation}%
and this relation determines uniquely $a_{i}(z)$ in terms of earlier
coefficients. Note that $a_{i}(z)$ is holomorphic in $D_{\kappa }$ and, by (%
\ref{a0knorm}) and (\ref{maj}) 
\begin{equation}
\sup_{z\in D_{\kappa }(0)}\left\vert a_{0}(z)\right\vert \leq \delta A~,
\label{delta}
\end{equation}%
by letting $\kappa $ small enough, for any $\delta >0$. Now, to obtain an
estimate on the growth rate of $|a_{i}(z)|$, let $\varphi _{i}$ denote the $%
i $-th Nagumo norm%
\begin{equation}
\Vert {}a_{i}\Vert _{i}:=\sup_{z\in {}D_{\kappa }(0)}\left( d_{\kappa
}(z)\right) ^{i}\left\Vert a_{i}(z)\right\Vert ,\hspace{6mm}\text{where}%
\hspace{3mm}d_{\kappa }(z)=\kappa -|z|.  \label{varphii}
\end{equation}%
of $a_{i}(z)$ and let $\beta _{n,m}$ the supremum in $D_{\kappa }$ of $%
\left\Vert B_{n,m}(z)\right\Vert $. The properties we shall use on Nagumo's
norms is proved in (\cite{BK}) and references therein and are sumarized by

\begin{enumerate}
\item $\Vert {}f+g\Vert _{k}\leq \Vert {}f\Vert _{k}+\Vert {}g\Vert _{k}$;

\item $\Vert {}fg\Vert _{k+l}\leq \Vert {}f\Vert _{k}\,\Vert {}g\Vert _{l}$;

\item $\Vert {}f^{\prime }\Vert _{k+1}\leq {}e(k+1)\Vert {}f\Vert _{k}$;

\item $\Vert {}f\Vert _{k}\leq \kappa \Vert {}f\Vert _{k-1}$,
\end{enumerate}

\noindent for any two functions $f$ and $g$ holomorphic in $D_{\kappa }$ and
nonnegative integers $k,l$.

\bigskip Let us assume that%
\begin{equation}
\varphi _{l}\leq \delta {}C_{l}\frac{1}{\nu ^{l}}  \label{f_l}
\end{equation}%
holds for $l=1,2,\ldots ,i-1$with $C_{l}=Al!/l^{2}$, for some positive
constants $\delta $ and $\nu $ to be determined. Similarly to (\ref{alphab})
and (\ref{alpha_nm}), 
\begin{equation}
\beta _{n,m}=\left\Vert B\right\Vert _{0}\leq \frac{C_{1}}{\rho _{1}^{n}\rho
^{m}}\leq \frac{\delta (1-b)}{c}\frac{\delta C_{n}}{\rho ^{n+m}}
\label{betanm}
\end{equation}%
holds for some $C_{1}<\infty $ and $\rho _{1}$ large enough. Then, it
follows by (\ref{ai}), (\ref{betanm}), (\ref{convk}) and the properties of
Nagumo norms%
\begin{eqnarray}
{}\varphi _{i} &\leq &\frac{c}{1-b}\left( \Vert {}z\Vert _{0}\,\Vert
{}a_{i-1}^{\prime }\Vert _{i}+\sum_{n=1}^{i}\sum_{m=1}^{\infty }\beta
_{n,m}\kappa ^{n}\sum_{\substack{ i_{1},\ldots ,i_{m}\geq 0:  \\ %
i_{1}+\cdots +i_{m}=i-n}}\varphi _{i_{1}}\cdots \varphi _{i_{m}}\right) 
\notag \\
&\leq &\frac{c}{1-b}{}e\kappa i{}\varphi _{i-1}+\delta \frac{1}{\nu ^{i}}%
\sum_{n=1}^{i}\sum_{m=1}^{\infty }C_{i-n}\left( \frac{\kappa }{\rho }\right)
^{n}C_{n}\left( \frac{\delta }{\rho }\right) ^{m}  \notag \\
&\leq &\left( \frac{2c}{1-b}e\kappa {}\nu +\frac{\delta }{\rho -\delta }%
\right) {}\delta C_{i}\frac{1}{\nu ^{i}}\leq \delta {}C_{i}\frac{1}{\nu ^{i}}%
,  \label{induction}
\end{eqnarray}%
where the last inequality holds provided $\delta <\rho /2$ and%
\begin{equation}
{}\nu \leq \frac{1-b}{2ce\kappa }\left( 1-\frac{\delta }{\rho -\delta }%
\right)  \label{one}
\end{equation}%
and this completes the induction: 
\begin{equation}
\sup_{z\in {}D_{\kappa }(0)}\left( d_{\kappa }(z)\right) ^{l}\left\vert
a_{l}(z)\right\vert \equiv \Vert {}a_{l}\Vert _{l}\leq \delta \frac{Al!}{%
l^{2}}\frac{1}{\nu ^{l}}\hspace{10mm}\forall \,l\geq 1.  \label{nninit}
\end{equation}%
with $\delta $ and $\nu $ fixed so that (\ref{delta}) and (\ref{one}) hold.

By definition (\ref{varphii}) of Nagumo norm,%
\begin{equation*}
\left\Vert a_{i}(z)\right\Vert \leq \frac{1}{(\kappa -\sigma )^{i}}\Vert
{}a_{i}\Vert _{i}\leq {}Ci!\mu ^{i}
\end{equation*}%
holds for all $i\geq 1$ uniformly in $\bar{D}_{\sigma }(0)$ for some $\sigma
<\kappa $, with $C=\delta {}A$ and $\mu ^{-1}=\nu (\kappa -\sigma )$, which
concludes the proof of Proposition \ref{P2}.

\hfill $\Box $

\section{Gevrey asymptotics\label{GA}}

\setcounter{equation}{0} \setcounter{theorem}{0}

In order to set up an equation involving derivatives of $f$ with respect to $%
\varepsilon $, we write

\begin{eqnarray*}
\phi _{i}(\varepsilon ,z) &=&\frac{1}{i!}\frac{\partial ^{i}f}{\partial {%
\varepsilon }^{i}}(\varepsilon ,z) \\
\phi _{i}^{\prime }(\varepsilon ,z) &=&\frac{\partial \phi _{i}}{\partial {z}%
}(\varepsilon ,z)~
\end{eqnarray*}%
and $\phi (\varepsilon ,z)=\left( \phi _{i}(\varepsilon ,z)\right) _{i\geq
0} $ for the sequence of those functions defined on ${}S(0,\gamma ;E)\times
D_{\kappa }(0)$; analogously to (\ref{Fab}) and (\ref{FB}), we write 
\begin{eqnarray*}
F(\varepsilon ,z,f) &=&\sum_{m=0}^{\infty }C_{m}(\varepsilon ,z)f^{m} \\
F^{[i,0,0]}(\varepsilon ,z,f) &=&\sum_{m=0}^{\infty
}C_{m}^{[i,0]}(\varepsilon ,z)f^{m}
\end{eqnarray*}%
for the $i$--th derivative of $h$ with respect to the first argument divided
by $i!$. The $i$-th total derivative of $F$ with respect to $\varepsilon $
can thus be written as 
\begin{eqnarray}
G_{i}(\varepsilon ,z,\phi _{0},\ldots ,\phi _{i}) &=&\frac{1}{i!}\frac{d^{i}%
}{d\varepsilon ^{i}}F(\varepsilon ,z,f)  \notag \\
&=&\sum_{m}\left( C_{m}^{[\cdot ,0]}(\varepsilon ,z)\ast \underset{m}{%
\underbrace{\phi (\varepsilon ,z)\ast \cdots \ast \phi (\varepsilon ,z)}}%
\right) _{i}  \notag \\
&=&T(\varepsilon ,z)\phi _{i}+\tilde{G}_{i}(\varepsilon ,z,\phi _{0},\ldots
,\phi _{i-1})  \label{TGi}
\end{eqnarray}%
where%
\begin{equation}
T(\varepsilon ,z)=\sum_{n=0}^{\infty }\sum_{m=1}^{\infty
}mA_{n,m}(\varepsilon )z^{n}\phi _{0}(\varepsilon ,z)^{m-1}  \label{T}
\end{equation}%
and $\tilde{G}_{i}(\varepsilon ,z,\phi _{0},\ldots ,\phi _{i-1})$ depends
only on derivatives of $f$ with respect to $\varepsilon $ of order lower
than $i$.

Differentiating equation (\ref{F}) $i$ times with respect to $\varepsilon $,
dividing by $i!$, we have 
\begin{equation}
\varepsilon z\phi _{i}^{\prime }-T(\varepsilon ,z)\phi
_{i}=H_{i}(\varepsilon ,z)  \label{THi}
\end{equation}%
for $i\geq 1$, where%
\begin{equation}
H_{i}(\varepsilon ,z)=\tilde{G}_{i}(\varepsilon ,z,\phi _{0},\ldots ,\phi
_{i-1})-z\phi _{i-1}^{\prime }~,  \label{Hi}
\end{equation}%
may be think as inhomogeneous holomorphic function of $\left( \varepsilon
,z\right) $ in $S(0,\gamma ;E)\times D_{\sigma }(0)$, and for $i=0$ simply (%
\ref{F}):%
\begin{equation*}
\varepsilon z\phi _{0}^{\prime }=F(\varepsilon ,z,\phi _{0})~.
\end{equation*}

\begin{proposition}
\label{P3} Let $f(\varepsilon ,z)$ be the unique holomorphic solution of (%
\ref{F}) on $S(0,\gamma ;E)\times D_{\sigma }(0)$ with $\sigma $, $\gamma $
and $E$ as in Proposition \ref{P1}. There exist $0<\sigma _{1}\leq \sigma $, 
$0<E_{1}\leq E$ and positive constants $C$ and $\mu $ such that%
\begin{equation*}
\left\Vert \phi _{i}(\varepsilon ,z)\right\Vert \leq {}Ci!\mu ^{i}
\end{equation*}%
holds for all $i\geq 0$ and every point $(\varepsilon ,z)$ in $S(0,\gamma
;E_{1})\times \bar{D}_{\sigma _{1}}(0)$.
\end{proposition}

\noindent \textbf{Proof} The case $i=0$ follows straightforwardly from
Proposition \ref{P1}. (\ref{THi}) is a linear singular perturbation equation
with regular singularity which can be dealt with the following auxiliary
result due to Balser-Kostov \cite{BK} (see Lemma 3 therein). For this, we
drop temporarily all subindices $i$ in (\ref{THi}).

Let%
\begin{equation}
T(\varepsilon ,z)-t_{0}(\varepsilon )=\sum_{n=1}^{\infty }{t_{n}(\varepsilon
)\,z^{n}}=S(\varepsilon ,z)  \label{AA}
\end{equation}%
and consider a sequence $\left( \psi _{k}(\varepsilon ,z)\right) _{k\geq 0}$
satisfying the system%
\begin{equation}
\left\{ 
\begin{array}{l}
\varepsilon {}z\psi _{0}^{\prime }(\varepsilon ,z)-t_{0}(\varepsilon )\psi
_{0}(\varepsilon ,z)=H(\varepsilon ,z)\vspace{2mm} \\ 
\varepsilon {}z\psi _{k}^{\prime }(\varepsilon ,z)-t_{0}(\varepsilon )\psi
_{k}(\varepsilon ,z)=S(\varepsilon ,z)\psi _{k-1}(\varepsilon ,z)\ ,\hspace{%
10mm}k=1,2,\ldots%
\end{array}%
\right. .  \label{system}
\end{equation}

By (\ref{AA}) and linearity, the sum over all equations in (\ref{system})
yields an equation of the form (\ref{THi}) satisfying by the sum $\psi
(\varepsilon ,z)=\displaystyle\sum_{k=0}^{\infty }\psi _{k}(\varepsilon ,z)$%
. We assume that $H(\varepsilon ,z)$ admits an expansion%
\begin{equation}
H(\varepsilon ,z)=\sum_{n=0}^{\infty }{h_{n}(\varepsilon )\,z^{n}\ }
\label{gn}
\end{equation}%
absolutely convergent for $\left\vert z\right\vert \leq \sigma $, uniformly
in $S(0,\gamma ;E)$. For $H$ given by (\ref{TGi}) and (\ref{Hi}) this will
actually be proven by induction when we resume the proof of Proposition \ref%
{P3}. We write, in addition, $f(z)\ll F(z)$ if $f(z)=\sum_{k=0}^{\infty }{%
c_{k}z^{k}}$ is majorized by $F(z)=\sum_{k=0}^{\infty }{C_{k}z^{k}}$ , i.
e., if $\left\vert c_{k}\right\vert \leq {}C_{k}$ holds for all $k$. If $f$
is a $\nu $-vector or a $\nu \times \nu $ matrix $f(z)\ll F(z)$ means
majorized relation for each component.

\begin{lemma}
\label{L4} There exist unique functions $\left( \psi _{k}(\varepsilon
,z)\right) _{k\geq 0}$, holomorphic in $S(0,\gamma ;E_{1})\times \bar{D}%
_{\sigma }(0)$, satisfying (\ref{system}). Each $\psi _{k}(\varepsilon ,z)$
has a zero of order $k$ at $z=0$: $\psi _{k}^{(0,k)}(\varepsilon ,0)=0$, and
satisfies%
\begin{equation}
\psi (\varepsilon ,z)=\sum_{k=0}^{\infty }{\psi _{k}(\varepsilon ,z)}\ll 
\frac{a}{I-a\Gamma (z)}\,\Omega (z)  \label{psi_major}
\end{equation}%
where%
\begin{eqnarray}
\Omega (z) &=&\sum_{n=0}^{\infty }\sup_{\varepsilon \in {}S(0,\gamma
;E)}\left\vert h_{n}(\varepsilon )\right\vert z^{n}~  \label{G} \\
\Gamma (z) &=&\sum_{n=0}^{\infty }\sup_{\varepsilon \in {}S(0,\gamma
;E)}\left\vert t_{n}(\varepsilon )\right\vert z^{n}~.  \label{H}
\end{eqnarray}%
holds for some $a<\infty $ provided $\sigma _{1}$ is small enough. $\psi
(\varepsilon ,z)$ is, in addition, the unique analytic solution in $%
S(0,\gamma ;E_{1})\times {}D_{\sigma }(0)$ of%
\begin{equation}
\varepsilon {}z\psi ^{\prime }(\varepsilon ,z)-T(\varepsilon ,z)\psi
(\varepsilon ,z)=H(\varepsilon ,z)  \label{psi_eq}
\end{equation}%
\vspace{1mm}with $\psi (\varepsilon ,0)=0$.
\end{lemma}

\noindent \textbf{Proof} Plugging%
\begin{equation*}
\psi _{k}(\varepsilon ,z)=\sum_{n=k}^{\infty }\varpi _{n,k}(\varepsilon
)\,z^{n}
\end{equation*}%
into (\ref{system}), yields%
\begin{eqnarray*}
(\varepsilon {}nI+t_{0}(\varepsilon ))\varpi _{n,0}(\varepsilon )
&=&h_{n}(\varepsilon )\ ,\qquad n\geq 0 \\
(\varepsilon {}nI+t_{0}(\varepsilon ))\varpi _{n,k}(\varepsilon )
&=&\sum_{m=k-1}^{n-1}t_{n-m}(\varepsilon )\,\varpi _{m,k-1}(\varepsilon )\ ,
\end{eqnarray*}%
for $1\leq k\leq n$~and $n\geq 1$. Observe that, by (\ref{T}) and (\ref{AA}%
), together with the fact that $\phi _{0}(\varepsilon ,0)=\displaystyle%
\sum_{j\geq 1}a_{j}(0)\varepsilon ^{j}$ (recall $a_{0}(0)=0$), 
\begin{equation*}
\varepsilon {}nI+t_{0}(\varepsilon )=\left( \varepsilon
{}nI+A_{0,1}(\varepsilon )\right) \left( I+\left( \varepsilon
{}nI+A_{0,1}(\varepsilon )\right) ^{-1}\sum_{m\geq 1}mA_{0,m}(\varepsilon
)\phi _{0}(\varepsilon ,0)^{m-1}\right)
\end{equation*}%
is invertible for every $\varepsilon \in D_{E_{1}}(0)$ if we take $E_{1}\leq
E$ so small that 
\begin{equation*}
c\sup_{\varepsilon \in \bar{D}_{E_{1}}(0)}\sum_{m=2}^{\infty }m\left\Vert
A_{0,m}(\varepsilon )\right\Vert \left\Vert \phi _{0}(\varepsilon
,0)\right\Vert ^{m-1}\leq d<1
\end{equation*}%
and $\left\Vert \left( \varepsilon {}nI+t_{0}(\varepsilon )\right)
^{-1}\right\Vert \leq c/(1-d)\equiv a<\infty $ holds uniformly in $%
D_{E_{1}}(0)$.

From these relations, we have%
\begin{equation*}
\psi _{0}(\varepsilon ,z)=\sum_{n=0}^{\infty }{\frac{1}{\varepsilon
{}nI+t_{0}(\varepsilon )}h_{n}(\varepsilon )\,z^{n}},
\end{equation*}%
and%
\begin{eqnarray*}
\psi _{k}(\varepsilon ,z) &=&\sum_{n=k}^{\infty }\frac{1}{\varepsilon
{}nI+t_{0}(\varepsilon )}\sum_{m=k-1}^{n-1}t_{n-m}(\varepsilon )\,\varpi
_{m,k-1}(\varepsilon )\,z^{n} \\
&=&\sum_{m=k-1}^{n-1}\left( \sum_{l=1}^{\infty }\frac{1}{\varepsilon
(m+l)I+t_{0}(\varepsilon )}t_{l}(\varepsilon )\,z^{l}\right) \varpi
_{m,k-1}(\varepsilon )\,z^{m}~.
\end{eqnarray*}%
Defining%
\begin{equation*}
\Psi _{k}(z)=\sum_{n=k}^{\infty }\sup_{\varepsilon \in {}S(0,\gamma
;E_{1})}\left\vert \varpi _{n,k}(\varepsilon )\right\vert \,z^{n},
\end{equation*}%
it follows, by (\ref{c}), (\ref{G}) and (\ref{H}) that%
\begin{eqnarray*}
\Psi _{0}(|z|) &\leq &{}a\Omega (|z|) \\
\Psi _{k}(|z|) &\leq &{}a\Gamma (|z|)\Psi _{k-1}(|z|)
\end{eqnarray*}%
for $k\geq 1$. Since $\psi _{k}(\varepsilon ,z)\ll \Psi _{k}(z)$ for $k\geq
1 $ and $\psi _{0}(\varepsilon ,z)\ll a\Omega (z)$ for $k=0$ hold for all $%
(\varepsilon ,z)\in {}S(0,\gamma ;E_{1})\times \bar{D}_{\sigma }(0)$, we
conclude (\ref{psi_major}) provided the geometric series $~\sum_{k\geq
1}a^{k}\left\Vert \Gamma (\sigma _{1})\right\Vert ^{k}$ converges. By (\ref%
{AA}) and (\ref{T})%
\begin{equation*}
\left\Vert \Gamma (\sigma _{1})\right\Vert =\sum_{n=1}^{\infty }{%
\sup_{\varepsilon \in {}S(0,\gamma ;E_{1})}}\left\Vert {t_{n-1}(\varepsilon )%
}\right\Vert {\,\sigma ^{n}}<\frac{1}{a}
\end{equation*}%
if $\sigma $ is small enough and thence, $\displaystyle\sum_{k=0}^{\infty }{%
\psi _{k}(\varepsilon ,}${$z$}${)}=\psi (\varepsilon ,z)$ is a uniformly
convergent series of analytic functions in $S(0,\gamma ;E_{1})\times
{}D_{\sigma }(0)$ which solves (\ref{psi_eq}). Since no other solution of (%
\ref{psi_eq}), regular at $z=0$, exists, the proof of Lemma \ref{L4} is
concluded.

\hfill $\Box $

We continue the proof of Proposition \ref{P3}. It remains to show that the
series (\ref{gn}) is uniformly convergent in $S(0,\gamma ;E_{1})\times
{}D_{\sigma }(0)$. This follows by induction. Clearly, $h_{0}(\varepsilon
,z) $ is holomorphic in $S(0,\gamma ;E_{1})\times {}D_{\sigma }(0)$. Suppose
that $\phi _{j}(\varepsilon ,z)$ is holomorphic in $S(0,\gamma ;E_{1})\times
{}D_{\sigma }(0)$ for each $1\leq j<i$. Then, by (\ref{Hi}), $%
h_{i}(\varepsilon ,z)$, is holomorphic in the same domain. By Lemma \ref{L4}%
, $\phi _{i}(\varepsilon ,z)$ is holomorphic in $S(0,\gamma ;E_{1})\times
{}D_{\sigma }(0)$ and, by (\ref{Hi}), we conclude it also holds for $%
h_{i+1}(\varepsilon ,z)$, justifying its representation as a convergent
series (\ref{gn}), uniformly in $S(0,\gamma ;E_{1})\times {}D_{\sigma }(0)$.
By induction, $\phi _{i}(\varepsilon ,z)$ is holomorphic in $S(0,\gamma
;E_{1})\times {}D_{\sigma }(0)$ for each $i\geq 1$ and%
\begin{equation}
\phi _{i}(\varepsilon ,z)\ll \frac{a}{I-a\Gamma (z)}\Omega _{i}(z)\ll \frac{a%
}{I-a\Gamma (\sigma _{1})}\Omega _{i}(z)\,.  \label{phi_major}
\end{equation}%
where $\Gamma _{i}$ depends on the $\phi _{j}(\varepsilon ,z)$ with $j<i$.
For $i=0$, by (\ref{G}),%
\begin{equation}
\left\vert f(\varepsilon ,z)\right\vert =\left\vert \phi _{0}(\varepsilon
,z)\right\vert \leq \frac{a}{I-a\Gamma (\sigma _{1})}\Omega _{0}(|z|)\leq
{}e_{0}  \label{e_0}
\end{equation}%
holds for all $\varepsilon \in {}S(0,\gamma ;E_{1})$ and $z\in {}D_{\sigma
}(0)$. For $i\geq 1$, we consider the modification of Nagumo norms:%
\begin{equation*}
\Vert {}f\Vert _{j}=\sup_{z\in {}D_{\sigma }(0)}\left( d_{\sigma
_{1}}(z)\right) ^{j}\sum_{n=0}^{\infty }\sup_{\varepsilon \in {}S(0,\gamma
;E)}\frac{1}{n!}\left\Vert \frac{\partial ^{n}f}{\partial {z}^{n}}%
(\varepsilon ,0)\right\Vert |z|^{n},
\end{equation*}%
with $d_{\sigma }(z)=\sigma -|z|$. It follows from (\ref{phi_major}) that%
\begin{equation}
\Vert \phi _{i}\Vert _{i}\leq \frac{a}{I-a\Gamma (\sigma _{1})}\Vert
H_{i}\Vert _{i},  \label{PhiGi}
\end{equation}%
where, by (\ref{TGi}),

\begin{equation*}
\Vert H_{1}\Vert _{1}\leq \sum_{m}\left\Vert C_{m}^{[1,0]}\right\Vert
_{1}\left\Vert \phi \right\Vert _{0}^{m}\leq C\frac{\left\Vert \phi
\right\Vert _{0}}{\rho -\left\Vert \phi \right\Vert _{0}}
\end{equation*}%
and, together with the properties of Nagumo norms, for $i\geq 2$ 
\begin{equation}
\Vert H_{i}\Vert _{i}\leq \Vert {}z\Vert _{0}\,\Vert \phi _{i-1}^{\prime
}\Vert _{i}+\sum_{m}\sum_{\substack{ i_{0},\ldots i_{m}\geq 0:  \\ %
i_{0}+\cdots +i_{m}=1}}\left\Vert C_{m}^{[i_{0},0]}\right\Vert
_{i_{0}}\left\Vert \phi _{i_{1}}\right\Vert _{i_{1}}\cdots \left\Vert \phi
_{i_{m}}\right\Vert _{i_{m}}~.  \label{Gi}
\end{equation}%
From these, together with (\ref{PhiGi}), a recursive relation of the same
type studied in Section \ref{FPS} may be derived for the $\Vert \phi
_{l}\Vert _{l}$ (see (\ref{f_l})-(\ref{nninit})) and one may conclude that%
\footnote{%
The details for this estimate are left to the reader.}%
\begin{equation}
\Vert \phi _{l}\Vert _{l-1}\leq \Delta \frac{Al!}{l^{2}}\omega ^{l}
\label{Phil}
\end{equation}%
holds for all $l\geq 1$ and some suitable constants $\Delta $ and $\omega $.
Picking $\sigma _{1}<\sigma $ together with the property of Nagumo norms,
yields 
\begin{equation*}
\left\vert \phi _{i}(\varepsilon ,z)\right\vert \leq \frac{\sigma }{(\sigma
-\sigma _{1})^{i}}\Vert {}\phi _{i}\Vert _{i-1}\leq {}Ci!\mu ^{i}
\end{equation*}%
for all $i\geq 1$ uniformly in $S(0,\gamma ;E_{1})\times \bar{D}_{\sigma
_{1}}(0)$, with $C=\sigma \Delta {}A$ and $\mu =\omega /(\sigma -\sigma
_{1}) $. We choose $C=\max (e_{0},\sigma _{1}\Delta {}A)$ in order to
include the $i=0$ case. This concludes the proof of Proposition \ref{P3}. 
\newline

\hfill $\Box $

\section{Summability\label{S}}

\setcounter{equation}{0} \setcounter{theorem}{0}

\begin{theorem}
\label{T1} Let (\ref{F}) be considered with $F$ given by (\ref{Fab}) where
the eigenvalues of $A_{0,1}(0)$ obey hypothesis (\ref{lambda}) for $%
(\varepsilon ,z)$ in a domain $S(0,\gamma ;E)\times D_{\sigma }(0)$ with $%
\gamma >\pi $. Then, there exist a radius $\sigma >0$ such that for $z\in 
\bar{D}_{\sigma }(0)$ the formal solution $\hat{f}(\varepsilon ,z)$ is $1$%
-summable in $\theta =0$ direction.
\end{theorem}

\noindent \textbf{Proof} By Taylor's Theorem%
\begin{equation*}
r_{I}(\varepsilon ,z)=\varepsilon ^{-I}\left( f(\varepsilon
,z)-\sum_{i=0}^{I-1}{f_{i}(z)\,\varepsilon ^{i}}\right) =\frac{I}{%
\varepsilon ^{I}}\int_{0}^{\varepsilon }{f_{I}(\zeta ,z)\,(\varepsilon
-\zeta )^{I-1}\,d\zeta },
\end{equation*}%
where the integral is along a path from $0$ to $\varepsilon $ inside $%
S(0,\gamma ;E)$. This, together with Proposition \ref{P3}, implies%
\begin{equation*}
\left\vert r_{I}(\varepsilon ,z)\right\vert \leq {}CI!^{s^{\prime }}\mu ^{I}
\end{equation*}%
for every $I$ and $(\varepsilon ,z)\in {}S^{\prime }\times \bar{D}_{\sigma
}(0)$, with $S^{\prime }$ any proper subsector of $S(0,\gamma ;E)$. In
addition, Proposition \ref{P2} states that $\hat{f}(\varepsilon ,z)$, a
formal solution of (\ref{F}), is an element of $\mathcal{O}(\sigma
)[[\varepsilon ]]_{1}$; therefore is an element of $\mathcal{O}(\sigma
)[[\varepsilon ]]_{1}$ for any $\sigma _{1}<\sigma $. Take now $\sigma _{1}$
and $E_{1}$ suffciently small. Hence, by definition (see Section 1.5 of \cite%
{Ba}), $\hat{f}(\varepsilon ,z)$ is an asymptotic expansion of order $1$, as 
$\varepsilon \rightarrow 0$ in the sector $S(0,\gamma ;E_{1})$, of $%
f(\varepsilon ,z)$, which by Proposition \ref{P1} is an analytic solution of
(\ref{F}) in the domain $S(0,\gamma ;E_{1})\times \bar{D}_{\sigma _{1}}(0)$.
Then, as $\gamma >\pi $, by hypothesis, $f(\varepsilon ,z)$ is the only
Gevrey order $1$ asymptotic expandable function in $S(0,\gamma ;E_{1})$
which has $\hat{f}(\varepsilon ,z)$ as its asymptotic expansion, and $\hat{f}%
(\varepsilon ,z)$ is $1$-summable in $\theta =0$ direction (see e.g. Section
3.2 of \cite{Ba}).

\hfill $\Box $

\begin{center}
\textbf{Acknowledgments}
\end{center}

DHUM thanks Gordon Slade and David Brydges for his hospitality at the
University of British Columbia.

\end{document}